\documentclass{amsart}

\usepackage{setspace}
\usepackage{asymptote}
\usepackage{mathrsfs}
\usepackage{array}
\usepackage{graphicx}
\usepackage{subfig}

\usepackage{manfnt} %%%%%%%%%%%%%%%%%%%%%% FIXME %%%%%%%%%%%%%%%%%%%%%%%%
\newcommand{\R}{\ensuremath{\mathbb{R}}}

\DeclareMathOperator{\sps}{ss}
\DeclareMathOperator{\pl}{pl}
\DeclareMathOperator{\s}{s}
\DeclareMathOperator{\crn}{cr}
\DeclareMathOperator{\br}{b}
\DeclareMathOperator{\tc}{tc}
\DeclareMathOperator{\sbr}{sb}

\newtheorem{lem}{Lemma}[section]
\newtheorem{thm}[lem]{Theorem}

\newtheorem{fact}[lem]{Fact}

\newtheorem{cor}[lem]{Corollary}
\theoremstyle{definition}
\newtheorem{Def}[lem]{Definition}

\newtheorem{rem}[lem]{Remark}

\begin{document}

\title{Planar and Spherical Stick Indices of Knots}
\date{\today}

\author[Adams]{Colin Adams}
\address{Colin Adams, Department of Mathematics and Statistics, Williams College, Williamstown, MA 01267}
\email{Colin.C.Adams@williams.edu}

\author[Collins]{Dan Collins}
\address{Dan Collins, Department of Mathematics, Princeton University, Princeton, NJ 08544}
\email{djc224@cornell.edu}

\author[Hawkins]{Katherine Hawkins}
\address{Katherine Hawkins, Department of Mathematics, Episcopal High School, P.O. Box 271299, Houston, TX 77277}
\email{khawkins@ehshouston.org}

\author[Sia]{Charmaine Sia}
\address{Charmaine Sia, Department of Mathematics, Harvard University, Cambridge, MA 02138}
\email{sia@mit.edu}

\author[Silversmith]{Rob Silversmith}
\address{Rob Silversmith, Department of Mathematics, 530 Church St., University of Michigan, Ann Arbor, MI 48109-1043}
\email{rsilvers@mich.edu}

\author[Tshishiku]{Bena Tshishiku}
\address{Bena Tshishiku, Department of Mathematics, University of Chicago, Chicago, IL 60637}
\email{tshishikub10@mail.wlu.edu}

\thanks{Support for this research was provided by NSF Grant DMS-0850577 and 
Williams College. The spherical stick index was first investigated by the
SMALL 2007 knot theory group, including Colin Adams, Nikhil Agarwal, William George, Rachel
Hudson, Trasan Khandhawit, Rebecca Winarski, and Mary Wootters. The authors
would like to thank Allison Henrich for her valuable input and assistance.}

\begin{abstract}
 The stick index of a knot is the least number of line segments required to build the knot in space. We define two analogous
 2-dimensional invariants, the planar stick index, which is the least number of line segments in the plane to build a projection,  and the spherical stick index, which is the least number of great circle arcs to build a projection on the sphere. We find bounds on these quantities in terms of other
 knot invariants, and give planar stick and spherical stick
 constructions for torus knots and for compositions of trefoils. In particular, 
 unlike most knot invariants,we show that the spherical stick index distinguishes between the granny and square knots, and that composing a nontrivial knot with a second nontrivial knot need not increase its spherical stick index.\end{abstract}

\maketitle

\section{Introduction}

The stick index $s[K]$ of a knot type $[K]$ is the smallest number of
straight line segments required to create a polygonal conformation of $[K]$ in
space. The stick index is generally difficult to compute. However, stick
indices of small crossing knots are known, and stick indices for
certain infinite categories of knots have been determined:
\begin{thm}[\cite{Jin}] 
  If $T_{p,q}$ is a $(p,q)$-torus knot with $p < q < 2p$, $s[T_{p,q}] = 2q$.
\end{thm}
\begin{thm}[\cite{AdamsStickNumber}] 
  If $nT$ is a composition of $n$ trefoils, $s[nT] = 2n+4$.
\end{thm} 

Despite the interest in stick index, two-dimensional analogues have not 
been studied in depth. In a recent paper, Adams and Shayler
\cite{AdamsShayler} defined a new invariant, the projective stick
index. We modify their definition slightly:

\begin{Def}
  A \emph{planar stick diagram} of a knot type $[K]$ is a closed polygonal
  curve in the plane, with crossing information assigned to self-intersections,
  that represents $[K]$.  The \emph{planar stick index} $\pl[K]$ of a knot type
  is the smallest number of edges in any planar stick diagram of $[K]$. 
\end{Def}

\begin{figure} 
  \subfloat[]{\label{fig:plspstrefoilA}
    \includegraphics[width=.25\textwidth]{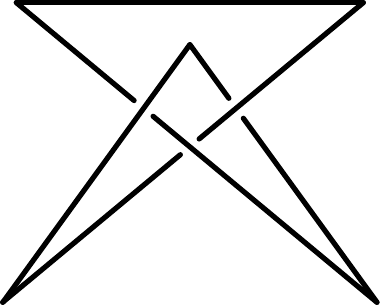}}
  \hspace{.2\textwidth}
  \subfloat[]{\label{fig:plspstrefoilB}
    \includegraphics[width=.25\textwidth]{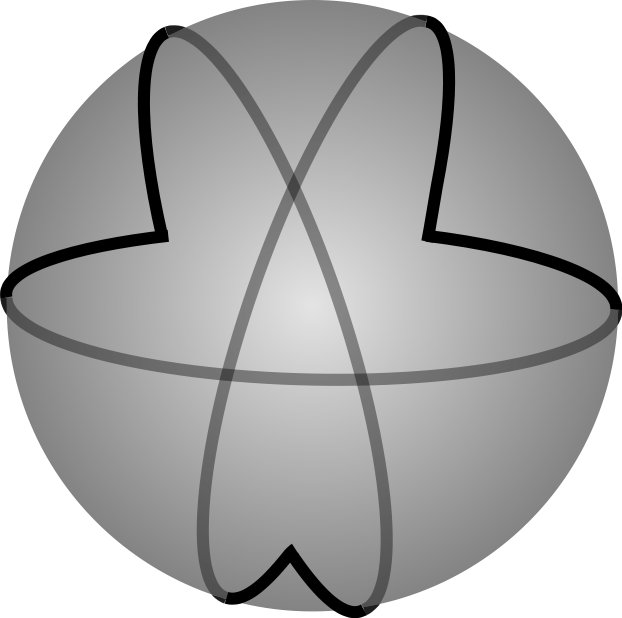}}
  \caption{A planar stick diagram and a spherical stick diagram of a trefoil.}
  \label{fig:plspstrefoil}
\end{figure}

An easy way to get a planar stick diagram of a knot type is to take a
3-dimensional stick conformation of that knot type and project it onto a plane.
Figure~\ref{fig:plspstrefoilA} shows a planar stick diagram of
a trefoil with five sticks.

We consider another invariant based on constructing diagrams of knots on the
sphere instead of in the plane.

\begin{Def}
  A \emph{spherical stick diagram} of $[K]$ is a closed curve on the sphere
  constructed from great circle arcs, with crossing information assigned to
  each self-intersection, that represents $[K]$.  The \emph{spherical stick
  index} $\sps[K]$ of a knot type is the minimum number of great circle arcs
  required to construct a spherical stick diagram of $[K]$.
\end{Def}

\begin{rem}
  We could define the spherical stick index of the unknot to be
  either $1$ or $2$, depending on whether we allow entire great
  circles in spherical stick diagrams. As such, we leave
  $\sps[\mbox{Unknot}]$ undefined. If we were to consider the
  spherical stick indices of links, the choice would become important.
\end{rem}

Figure~\ref{fig:plspstrefoilB} shows a spherical stick diagram of a trefoil.
A spherical stick diagram can be obtained via radial projection of a
stick knot in space onto a sphere from some point in space, or via
radial projection of a planar stick diagram from some point not in the
plane. 

In Section \ref{section:pl}, we establish bounds for the planar stick 
index in terms of other invariants, including crossing number, stick index,
and bridge index. Section \ref{section:sps} establishes similar bounds
for the spherical stick index.

In Section \ref{section:examples}, we construct planar stick diagrams
and spherical stick diagrams for torus knots and compositions of
trefoils, providing upper bounds for the planar stick index and
spherical stick index of these knot types. In some cases, the bounds
from Sections \ref{section:pl} and \ref{section:sps} show that the
constructions are minimal. 

Our results are as follows. Let $T_{p,q}$ denote the $(p,q)$-torus knot. 
\begin{thm} \label{thm:pltorus}
  Let $2\le p<q$. Then
  \begin{alignat*}{2}
   \pl[T_{p,q}] & \le 2q - 1 & \qquad & \textrm{if} \quad q < 2p ,\\
    \pl[T_{p,q}] & \le q & \qquad &\textrm{if} \quad 2p< q, \\
    \pl[T_{p,q}] & \ge 2p+1 & & \textrm{for all } p.
 \end{alignat*}
  When $q = p+1$ or $q = 2p+1$, the inequalities exactly determine the
  planar stick index to be $\pl[T_{p,p+1}] = \pl[T_{p,2p+1}] = 2p+1$.
\end{thm}
\begin{thm} \label{thm:spstorus}
  Let $2\le p < q$. Then
  \[ \sps[T_{p,q}] \le q. \]
  Moreover, $\sps[T_{q-1,q}]=q$.
\end{thm}

Let $nT$ denote a composition of $n$ trefoils (of any combination of
handedness), and $aT_L\# bT_R$ denote the composition of $a$ left-handed
trefoils with $b$ right-handed trefoils. Because composition of knots is
commutative and associative (see \cite{KnotBook}), $aT_L\# bT_R$ is
well-defined. 
\begin{thm} \label{thm:pltrefoil} 
  For $n \ge 1$, $ \pl[nT]=2n+3. $
\end{thm}
\begin{thm} \label{thm:spstrefoil}
  For $m \ge 1$, $\sps[mT_L\# mT_R] \le 2m+2. $
  \medskip
  
  \noindent More generally, for $0\le m < n$, 
  $\sps[mT_L\# nT_R] = \sps[nT_L\# mT_R] \le 2n+1. $
\end{thm}

The difference between Theorems \ref{thm:pltrefoil} and \ref{thm:spstrefoil} is
striking: the planar stick index of a composition of trefoils is independent
of the handedness of the trefoils composed, while our construction of a
spherical stick diagram depends heavily on handedness. It would be interesting
to know if the bounds in Theorem \ref{thm:spstrefoil} are sharp, and whether the
spherical stick index of a composition of trefoils depends on handedness in
general. This seems difficult to prove, since most invariants that we could use
to obtain lower bounds do not detect handedness of composites. However, by
classifying all knots with $\sps[K]\le 4$ (as we do in Section
\ref{section:sps4}), we can prove that the bound in Theorem
\ref{thm:spstrefoil} is sharp in the case of composing two trefoils.

\begin{thm} \label{thm:sps4classification}
  The nontrivial knot types with $\sps[K]\le 4$ are
  \[ 3_1,4_1,5_1,5_2,6_1,6_2,6_3,7_4,8_{18},8_{19},8_{20},\]
  and the square knot $T_L\# T_R$. All of these knots except the trefoil $3_1$
  have $\sps[K] = 4$.
\end{thm}
\begin{cor} \label{thm:spssquaregranny}
  $\sps[T_L\# T_R] = 4$, while $\sps[T_L\# T_L] = \sps[T_R\# T_R] = 5$.
\end{cor}

We also see a very unusual characteristic for a naturally defined physical knot invariant:

\begin{cor} \label{cor:spscomposition}
  There exist nontrivial knots $K_1$ and $K_2$ so that 
  \[ \sps[K_1\# K_2] = \sps[K_1]. \]
\end{cor}

%%%%%%%%%%%%%%%%%%%%%%%%%%%%%%%%%%%%%%%%%%%%%%%%%%%%%%%%%%%%%%%%%%%%%%%%%%%%%%
%%%%%%%%%%%%%%%%%%%%%%%%%%%%%%%%%%%%%%%%%%%%%%%%%%%%%%%%%%%%%%%%%%%%%%%%%%%%%%
\section{Planar Stick Index} \label{section:pl}
In general, the planar stick index of a knot is difficult to compute. It is 
straightforward to construct a planar stick diagram, but hard to prove that it
is minimal. In this section, we establish bounds on the planar stick index of
a knot in terms of other invariants. These bounds enable us to compute exact
values for planar stick index for certain categories of knots in Section \ref{section:examples}.

\begin{thm} \label{thm:pl-s}
  $\pl[K]\leq\s[K]-1$.
\end{thm}
\begin{proof}
  Consider a polygonal conformation of $[K]$ that realizes the stick index. If
  we project the knot onto a plane normal to one stick, that stick projects to a
  single point. In the ``generic case,'' the resulting polygonal curve in the
  plane is a diagram of $[K]$ with at most $s[K]-1$ edges.  The diagram
  fails to be generic if three edges intersect at the same point, or if a
  vertex overlaps an edge. In such a case, however, we can tweak the original
  conformation slightly so that after projecting, we obtain a generic
  $(s[K]-1)$-edge diagram of $[K]$. 
\end{proof}

\begin{thm} \label{thm:pl-crn}
  Let $\crn[K]$ be the crossing number of $[K]$. Then
	\[ \frac{3+\sqrt{9+8\crn[K]}}{2}\leq\pl[K]. \]
\end{thm}
\begin{proof}
  Consider a planar stick diagram of $[K]$ with $n=\pl[K]$ sticks.  Each stick
  can cross at most $n-3$ other sticks, since it can cross neither itself nor
  the two adjacent sticks. The total number of self-intersections (which is at
  least $\crn[K]$) is bounded above by $\frac{1}{2}n(n-3)$. Rearranging $\crn[K]
  \le \frac{1}{2}n(n-3)$ gives the desired inequality.
\end{proof}
\begin{thm} \label{thm:pl-comp}
  If $K_1\# K_2$ is the composition of two knots $K_1,K_2$, then
  \[ \pl[K_1\# K_2] \le \pl[K_1] + \pl[K_2] - 2. \]
\end{thm}
\begin{proof}
  Consider planar stick diagrams for $K_1$ and $K_2$ that realize planar stick
  index. Since any two adjacent sticks in the diagram of $K_1$ are
  non-parallel, we can perform an orientation-preserving linear transformation
  that makes two adjacent sticks of $K_1$ perpendicular. We do likewise for the
  diagram of $K_2$.
  
  Once we have these diagrams, we can rotate and attach them at the right
  angles such that the incident sticks line up.  The point at which the corners
  were attached then becomes a crossing. If $K_1$ and $K_2$ do not overlap, 
  the resultant diagram is a $(\pl[K_1]+\pl[K_2]-2)$-stick representation of
  $K_1\# K_2$.

  If $K_1$ and $K_2$ do overlap, first note that if necessary we can tweak the
  diagrams slightly so that the new diagram is generic. We can then choose the new
  crossings so that sticks from $K_1$ always cross over sticks from $K_2$. It
  is clear that the diagram represents $K_1\#K_2$.
\end{proof}

We get another bound on the planar stick index in terms of the
\emph{bridge index}. Let $\br(K,p)$ be the number of local maxima of
a knot conformation $K$ relative to a direction (taken to be a vector
$p$ on the 2-sphere $S^2$). The bridge index is given by
\[ \br[K] = \min_{K\in [K]} \min_{p\in S^2} \br(K,p). \]

This definition is similar to Milnor's definition of crookedness (see 
\cite{Milnor}). However, if an extremum occurs at an interval of constant
height, we count it as one extremum rather than infinitely many.

\begin{thm}\label{thm:pl-br}
  $2\br[K]+1 \le \pl[K]$
\end{thm}
\begin{proof}
  For a planar stick diagram, the total curvature is the sum of the exterior 
  angles. There are $\pl[K]$ vertices in a minimal planar stick diagram of a
  knot $[K]$. Since each vertex has an exterior angle strictly less than
  $\pi$, the total curvature of such a diagram is less than $\pi\pl[K]$.

  We view this diagram as a curve in a plane in 3-space.  Bending the
  sticks slightly out of the plane at each crossing yields a conformation of
  the knot (as opposed to a diagram). Since we can bend the sticks by an
  arbitrarily small amount, the final total curvature can be made arbitrarily
  close to the original total curvature. Since the original total curvature was
  strictly less than $\pi\pl[K]$, the final total curvature can be made to be
  less than $\pi\pl[K]$.

  Milnor showed in \cite{Milnor} that for any conformation $K$ with total
  curvature $\tc(K)$, 
  \[ 2\pi\br[K]<\tc(K). \] 
  Since $\tc(K)<\pi\pl[K]$, we find $2\br[K]<\pl[K]$.
  Since both quantities are integers, the result follows.
\end{proof}

%%%%%%%%%%%%%%%%%%%%%%%%%%%%%%%%%%%%%%%%%%%%%%%%%%%%%%%%%%%%%%%%%%%%%%%%%%%%%%
%%%%%%%%%%%%%%%%%%%%%%%%%%%%%%%%%%%%%%%%%%%%%%%%%%%%%%%%%%%%%%%%%%%%%%%%%%%%%%
\section{Spherical Stick Index} \label{section:sps}

\begin{figure}
  \includegraphics[height=.22\textwidth]{fig_spstrefoil.jpg}
  \hspace{.15\textwidth}
  \includegraphics[height=.22\textwidth]{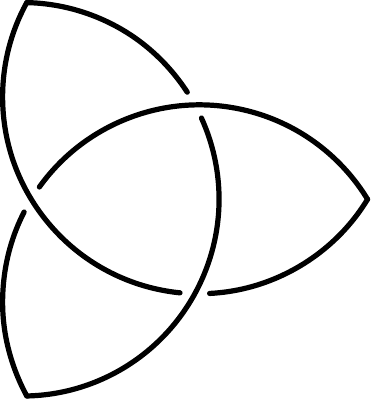}
  \caption{A spherical stick trefoil and its stereographic projection.}
  \label{fig:stereotrefoil}
\end{figure}

When studying spherical stick diagrams, it is helpful to consider their 
stereographic projections. Given a diagram of a knot on a sphere, 
we choose a point on the sphere not on the diagram to label $\infty$. Consider this the north pole. The
stereographic projection relative to this point maps $S^2\backslash\{\infty\}$ 
homeomorphically to $\R^2$, which is the plane though the equator, and transfers the diagram into $\R^2$. Moreover, 
stereographic projection preserves the knot type of a diagram. See Figure 
\ref{fig:stereotrefoil} for an example. The following fact can be proved relatively easily.

\begin{fact}  
  Stereographic projection gives a one-to-one correspondence between great 
  circles on the sphere that do not pass through infinity, and circles in the
  plane that have a diameter with endpoints $p,q$ that contains the origin, and
  satisfies $|p|\cdot|q| = 1.$
\end{fact} 

In some situations, it is more convenient to think about circles
in the plane rather than great circles on the sphere. Most of our figures
of spherical stick diagrams will be stereographic projections
for clarity. 

We prove bounds on spherical stick index, many of
which are analogous to those proven in Section \ref{section:pl} for planar
stick index.

\begin{thm} 
  $\sps[K]\le\pl[K]$.
\end{thm}
\begin{proof}
  Observe that given a polygonal curve in space, radial projection onto a
  sphere maps each edge to a great circle arc. 
  
  Consider a planar stick diagram that realizes planar stick index for $[K]$.
  We put it in space in a plane not containing the origin and radially project
  to the unit sphere. The diagram projects to a spherical stick diagram of
  $[K]$ with $\pl[K]$ great circle arcs.  
\end{proof}
\begin{thm} \label{thm:sps-s-pl}
  $\sps[K]\leq\s[K]-2$.
\end{thm}
\begin{proof}
  Consider a minimal stick realization of a knot in space. Using a similar trick as appears in \cite{Calvo},
  we choose a vertex
  $v$ of the knot, and radially project the knot (minus $v$) onto a sphere
  centered at $v$. Radial sticks project to points, and non-radial sticks
  project to great circle arcs. Since the two sticks adjacent to $v$ are
  radial, the projection has at most $s[K]-2$ arcs. However, it is no longer a
  closed curve, as there are two ``loose ends'' corresponding to the sticks
  incident at $v$.

  As projections of line segments, the great circle arcs must be strictly
  smaller than $\pi$ radians. Since any pair of distinct great circles
  intersect at two antipodal points, and each arc traverses less than half of a
  great circle, no two arcs can intersect more than once. In particular, the
  arcs with ``loose ends'' intersect at most once. We extend these arcs until
  they meet, making the extended arcs understrands at any newly-created
  crossings (to preserve the knot type). This yields a spherical stick diagram
  of the knot with $\s[K]-2$ arcs.
\end{proof}

\begin{thm} \label{thm:sps-crn}
  $1+\sqrt{1+\crn[K]}\leq\sps[K]$.
\end{thm}
\begin{proof}
  Consider a spherical stick diagram of $[K]$ with $n = \sps[K]$ great circle
  arcs. No arc intersects itself, and each arc intersects each of the other
  $n-1$ arcs at most twice. Also, vertices connecting arcs at endpoints are
  intersections, but are not crossings in the diagram. Thus, there are at most
  $2(n-1)-2 =2n-4$ crossings on each arc. Since each crossing occurs on exactly
  two arcs, there are at most
  \[ \frac{1}{2} n (2n-4) = n(n-2) \]
  crossings in the diagram, so $\crn[K]\le \sps[K](\sps[K]-2)$. Solving for
  $\sps[K]$ gives the desired inequality.
\end{proof}

\begin{figure} 
  \subfloat[]{\label{fig:sps-compA}
    \includegraphics[width=.3\textwidth]{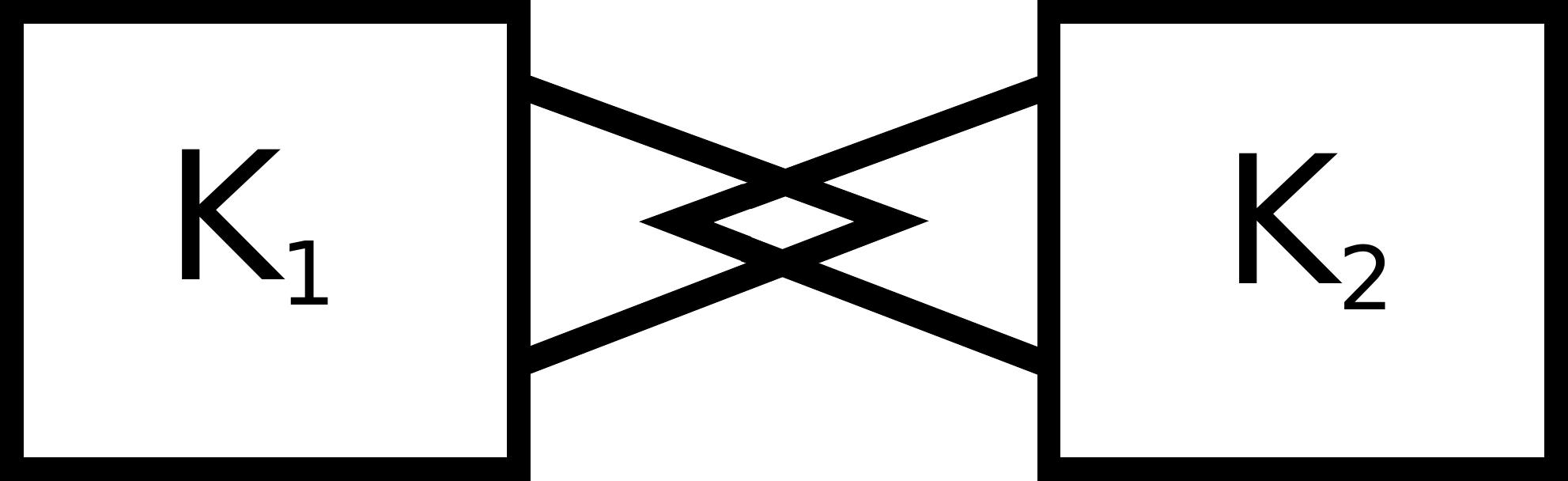}}
  \hspace{.15\textwidth}
  \subfloat[]{\label{fig:sps-compB}
    \includegraphics[width=.3\textwidth]{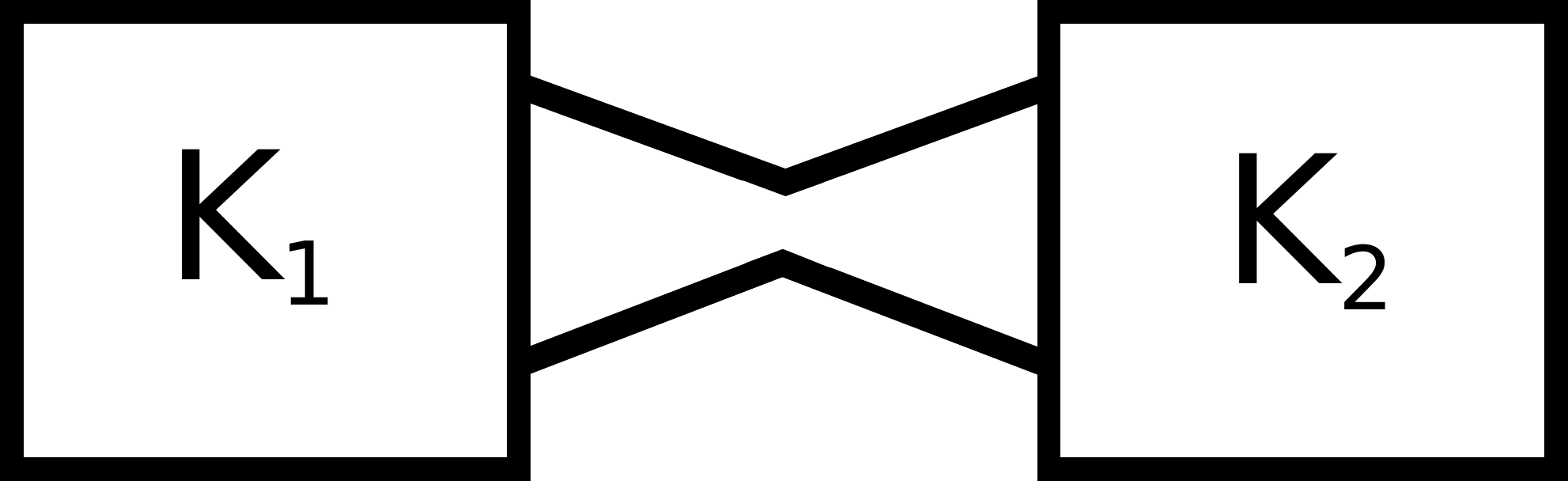}}
  \caption{Composition of spherical stick diagrams. Great circle arcs appear
    locally as straight lines.}
  \label{fig:sps-comp}
\end{figure}

\begin{thm} \label{thm:sps-comp}
  $\sps[K_1\# K_2] \le \sps[K_1]+\sps[K_2]$. 
\end{thm}
\begin{proof}
  Suppose we have minimal spherical stick diagrams of $K_1$ and $K_2$ on the
  sphere. We position them so that two vertices of the diagrams overlap as
  shown in Figure~\ref{fig:sps-compA}. Note that the diagrams
  may overlap in many other places. As before, we move $K_2$ to ensure that the
  diagram is generic, and choose the new crossings so that arcs in $K_1$ cross
  above arcs in $K_2$. We then change the diagrams as shown in Figure~\ref{fig:sps-compB} to obtain a diagram of $K_1\# K_2$. Our new
  diagram has $\sps[K_1]+\sps[K_2]$ great circle arcs.
\end{proof}

We can bound spherical stick index in terms of an invariant related to the
bridge index. Using $\br(K,p)$ as previously defined, we let the
\emph{superbridge index} of a knot $[K]$ be 
\[ \sbr[K] = \min_{K\in [K]} \max_{p\in S^2} \br(K,p), \] 
as in \cite{Kuiper}. The strict inequality $\br[K] < \sbr[K]$ holds for all
knot types, as proven in \cite{Kuiper}.

\begin{figure}
  \includegraphics[width=.2\textwidth]{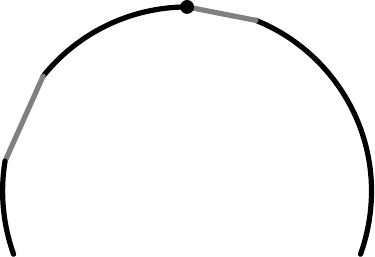}
  \hspace{.15\textwidth}
  \includegraphics[height=.2\textwidth]{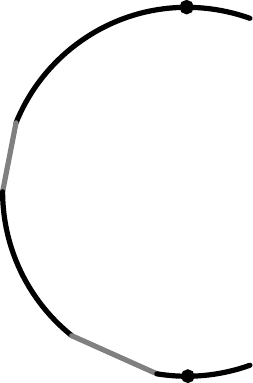}
  \caption{Almost-circular arcs from the proof of Theorem \ref{thm:sps-sb}} 
  \label{fig:sps-sb}
\end{figure}

\begin{thm} \label{thm:sps-sb}
  $\frac{2}{3}\sbr[K] + \frac{1}{3} \le \sps[K]$.
\end{thm}
\begin{proof}
  Consider a spherical stick diagram of $[K]$ with $n=\sps[K]$ great circle
  arcs.  We modify it as follows to obtain a conformation $K$ in space.
  For any crossing of the diagram, there is an ``overstrand'' and an
  ``understrand'', relative to the outside of the sphere. We remove a small
  portion of the understrand, replacing it with a straight line. We do this for
  all crossings, and obtain a conformation of $[K]$ that radially projects to
  our diagram. The conformation consists of $n$ almost-circular arcs, like
  those in Figure~\ref{fig:sps-sb}.

  Given a direction $p\in S^2$, we want an upper bound on the number of extrema 
  in the direction $p$. Each of the $n$ points connecting two almost-circular
  arcs can be an extremum. Other extrema must occur on the interiors of the
  arcs, and each arc can have at most two interior extrema (see Figure 
 ~\ref{fig:sps-sb}). Thus $K$ has at most $3n$ extrema in the 
  direction $p$. Since $\br(K,p)$ counts the number of maxima, 
  $\br(K,p)\le 3n/2$. Therefore,
  \[ \sbr[K] \le \max_{p\in S^2} \br(K,p) \le \frac{3n}{2} 
      = \frac{3}{2}\sps[K]. \]
  Rearranging gives $\frac{2}{3}\sbr[K]\le\sps[K]$. 

  To improve the bound by $1/3$, we use a small trick that guarantees an arc
  of length less than $\pi$. We stereographically project our original
  spherical stick diagram from a point not in the diagram whose antipode is in
  the diagram. We get a diagram in the plane consisting of circles and one
  line segment through the origin. We scale the diagram so the line segment is
  contained in the unit disc, and then stereographically project back to the 
  sphere. The result is a spherical diagram of $[K]$ with one great circle
  arc of length less than $\pi$, and $n-1$ other circular arcs (not 
  necessarily great circle arcs). We change these into almost-circular 
  arcs to obtain a conformation $K$, and apply the same counting argument as
  before. Because the great circle arc of length less than $\pi$ can have at
  most one interior extremum, we find
  \[ \sbr[K] \le \frac{3\sps[K]-1}{2}, \]
  which rearranges to the desired inequality. 
\end{proof}

Since bridge number is known for many more knots than is superbridge number, we note that Kuiper's result that $\br[K] < \sbr[K]$ implies:

\begin{cor}   $\frac{2}{3}\br[K] + 1 \le \sps[K]$.
\end{cor}

\begin{rem}
  In all cases where we know both the superbridge index and spherical
  stick index, $\sbr[K]\le\sps[K]$. It would be interesting to know
  whether this holds in general, as it would prove that some of the
  constructions in the next section realize spherical stick index. (See Remarks 4.2 and 4.4.)
\end{rem}

%%%%%%%%%%%%%%%%%%%%%%%%%%%%%%%%%%%%%%%%%%%%%%%%%%%%%%%%%%%%%%%%%%%%%%%%%%%%%%
%%%%%%%%%%%%%%%%%%%%%%%%%%%%%%%%%%%%%%%%%%%%%%%%%%%%%%%%%%%%%%%%%%%%%%%%%%%%%%
\section{Examples} \label{section:examples}
\subsection{Torus Knots}
Armed with bounds on planar and spherical stick indices, we examine some
classes of knots, and prove Theorems \ref{thm:pltorus}, \ref{thm:spstorus},
\ref{thm:pltrefoil}, and \ref{thm:spstrefoil} stated in the introduction. We
begin with the torus knots, one of the most easily described and exhaustively
studied classes of knots.

We need a few well-known properties of torus knots.  It is known that for any
$p$ and $q$, the $(p,q)$-torus knot is equivalent to the $(q,p)$-torus knot
(see \cite{KnotBook}). So, without loss of generality, we always assume $p<q$.
Also, we require $p$ and $q$ to be coprime (otherwise we get a torus link with
$\gcd(p,q)$ components).

We need the values of some invariants of $T_{p,q}$ (for $p < q$). First,
the bridge index has been shown (see \cite{Schubert} or \cite{Kuiper}) to be 
\[ \br[T_{p,q}] = p. \]
It was shown in \cite{Murasugi} that the crossing number is 
\[ \crn[T_{p,q}]=(p-1)q. \]
Finally, it was proven in \cite{Jin} that, for $q < 2p$, 
\[ \s[T_{p,q}] = 2q. \]

We prove Theorems \ref{thm:pltorus} and \ref{thm:spstorus}, which pertain
to planar and spherical stick indices of torus knots. 

\begin{figure}
  \subfloat[]{\label{fig:pltorusA}
    \includegraphics[width=.3\textwidth]{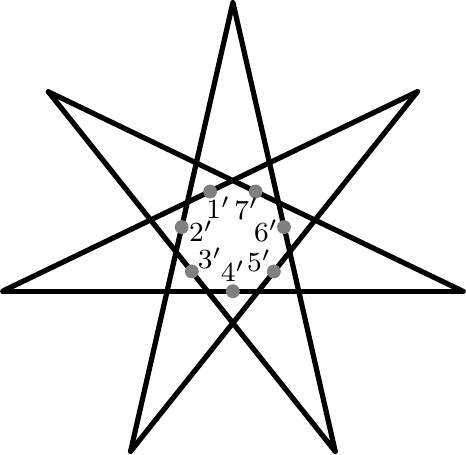}}
  \hspace{.01\textwidth}
  \subfloat[]{\label{fig:pltorusB}
    \includegraphics[width=.3\textwidth]{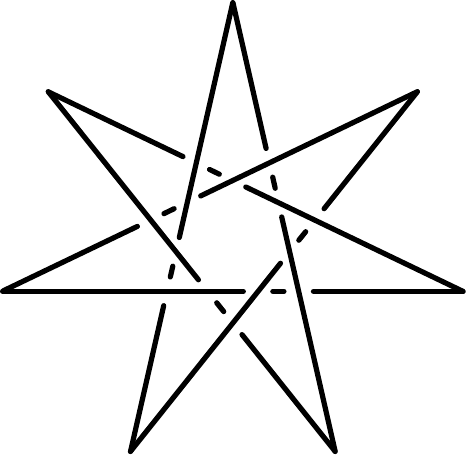}}
  \hspace{.01\textwidth}
  \subfloat[]{\label{fig:pltorusC}
    \includegraphics[width=.3\textwidth]{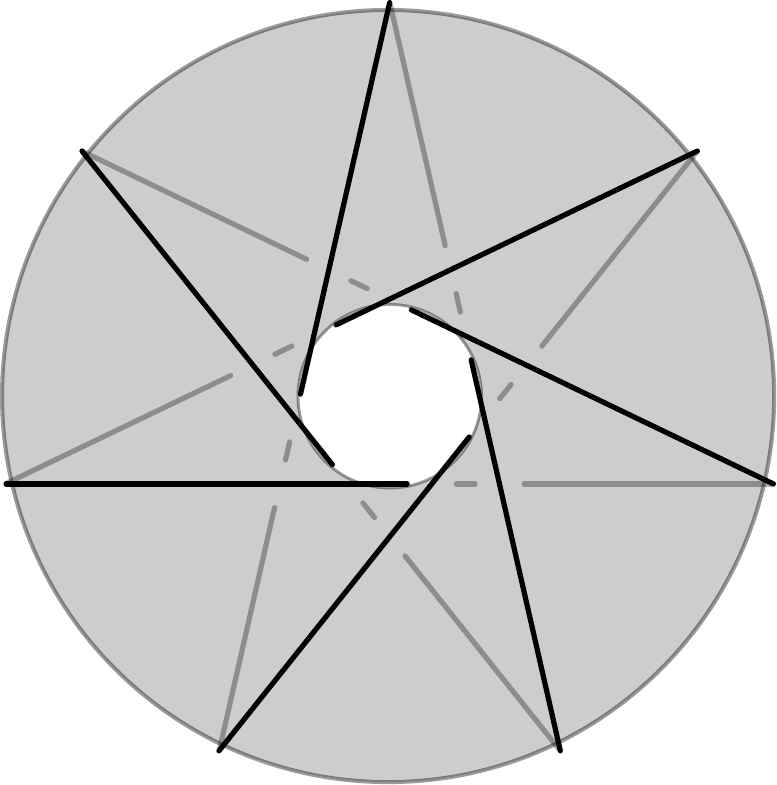}}
  \caption{Construction of $T_{3,7}$ with 7 sticks.}
  \label{fig:pltorus}
\end{figure}

\begin{proof}[Proof of Theorem \ref{thm:pltorus}]
  Two of the inequalities follow directly from theorems established in
  Section \ref{section:pl} and the facts above. We can apply Theorem
  \ref{thm:pl-s} to show that when $q < 2p$, 
  \[ \pl[T_{p,q}] \le 2\s[K]-1 = 2q-1. \]
  Similarly, Theorem \ref{thm:pl-br} gives
 \[ \pl[T_{p,q}] \ge  2\br[T_{p,q}] + 1 = 2p+1. \]
  
  It remains to show that for $2p < q$ we can construct a planar stick diagram
  of $T_{p,q}$ with $q$ sticks. We consider $q$ evenly spaced points on a
  circle, $z_1\ldots z_q$, labeled counterclockwise.  We then draw $q$ line
  segments, connecting $z_n$ to $z_{n+p}$ for each $n$. The result is a
  $q$-pointed star, as in Figure~\ref{fig:pltorusA}. We label the
  stick from $z_n$ to $z_{n+p}$ as stick $n$ (and take these labels modulo
  $q$). We label the midpoint of stick $n$ as $n'$. We can see that the middle
  of the diagram is a regular $q$-gon for which the midpoint of each side is
  some $n'$.
  
  By construction, stick $n$ attaches to stick $n+p$, which attaches to stick
  $n+2p$, and so on. Since $p$ and $q$ are coprime for torus knots, the $q$
  sticks form a single closed curve.  Furthermore, stick $n$ and stick $m$
  cross if and only if $m$ is in the following set (modulo $q$)
  \[ \{ n-p-1,\ldots,n-1,n+1,\ldots,n+p-1 \}. \]
  In this case, we let $I_{m,n} = I_{n,m}$ denote the point of intersection.
  For each stick, we define directions of ``clockwise'' and
  ``counterclockwise'' relative to the origin. On stick $n$, the intersections
  $I_{n,n-1},\ldots,I_{n,n-p+1}$ are clockwise from $n'$, and
  $I_{n,n+1},\ldots,I_{n,n+p-1}$ are counterclockwise from $n'$ (see Figure
 ~\ref{fig:pltorus}). This means our diagram has $(p-1)q$ intersections, which
  is equal to the crossing number of the standard projection of $T_{p,q}$. 

  We specify the crossings by letting stick $n$ be the overstrand for the $p-1$
  crossings clockwise from $n'$, and the understrand for the $p-1$ crossings
  counterclockwise from $n'$ (see Figure~\ref{fig:pltorusB}). 
  
  From Figure~\ref{fig:pltorusC}, we can see our diagram is a
  projection of a knot on the standard torus in $\R^3$. This knot winds around
  the torus $p$ times in one direction and $q$ times in the other, so it is
  $T_{p,q}$. 
\end{proof} 

\begin{rem}
  The fact that $\pl[T_{p,p+1}] = \pl[T_{p,2p+1}] = 2p+1$ is 
  striking. We suspect, but have not been able to prove, that 
  $\pl[T_{p,q}] = 2p+1$ whenever $p+1\le q\le 2p+1$. 
\end{rem}

\begin{figure}
  \includegraphics[width=.32\textwidth]{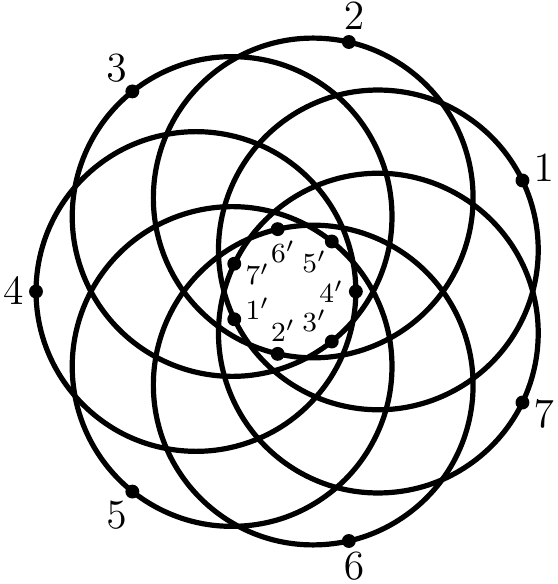}
  \caption{Seven great circles evenly spaced, for construction of 7-stick 
    spherical projective stick diagrams of $T_{p,7}$}
  \label{fig:spstorus_circles}
\end{figure}

\begin{proof}[Proof of Theorem \ref{thm:spstorus}]
  To show $\sps[T_{p,q}]\le q$, we use a similar construction as we used in the
  previous proof. We begin with a regular $q$-gon on the sphere, centered at
  the north pole. We extend the sides to great circles.  The stereographic
  projection is shown for $q=7$ in Figure~\ref{fig:spstorus_circles}. We label
  the circles counterclockwise from $C_1$ to $C_q$. We let the
  \emph{basepoint} of each circle be the farthest point on the circle from the
  origin (under the stereographic projection) and label it $1$ through $q$ as in the figure,
   and label the point antipodal to
  the basepoint $n$ as $n'$. As before, we consider these labels
  modulo $q$. 

  Next, we establish some facts about intersections of the circles. Any two
  great circles in the sphere intersect twice, at antipodal points. (The
  intersections are no longer antipodal under stereographic projection.) If we
  fix a circle $C_n$ in our diagram, any other circle $C_m$ intersects it once
  clockwise from $n$ and once counterclockwise from $n$. Furthermore, an
  intersection point of circles $C_m$ and $C_n$ that is clockwise from $n$ must be
  counterclockwise from $m$.  Hence for any $m\ne n$, there is a unique
  intersection $i_{m,n}$ of circle $C_m$ and circle $C_n$ so that $i_{m,n}$ is
  counterclockwise from $m$ and clockwise from $n$. 

\begin{figure}
  \subfloat[]{\label{fig:spstorusA}
    \includegraphics[width=.25\textwidth]{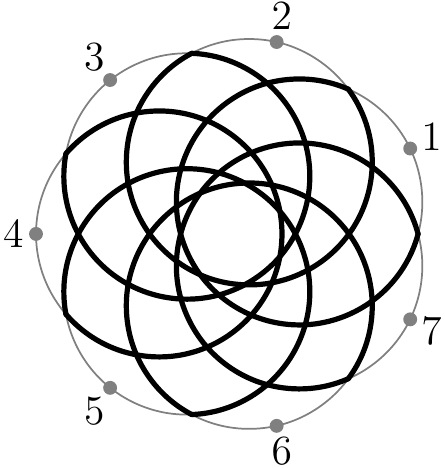}}
  \hspace{.05\textwidth}
  \subfloat[]{\label{fig:spstorusB}
    \includegraphics[width=.25\textwidth]{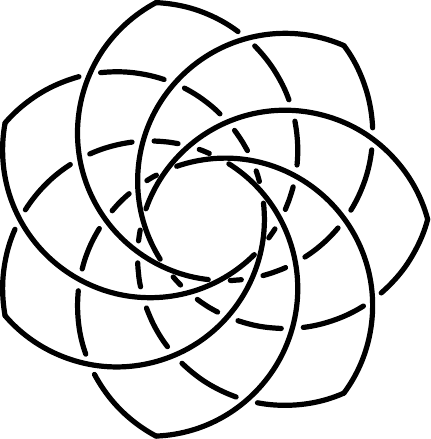}}
  \hspace{.05\textwidth}
  \subfloat[]{\label{fig:spstorusC}
    \includegraphics[width=.25\textwidth]{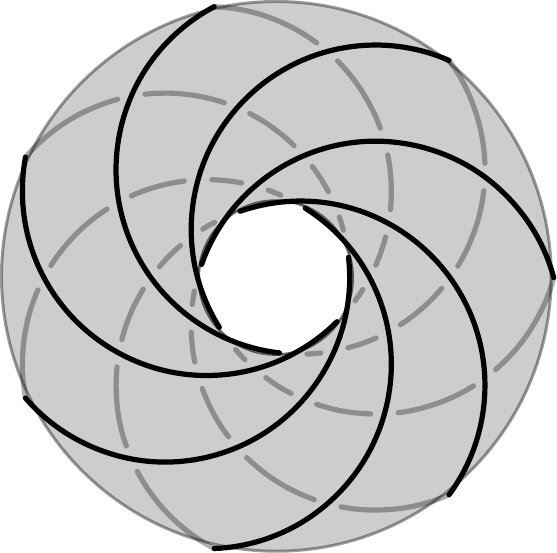}}
  \caption{Construction of $T_{6,7}$ with $7$ great circle arcs.}
  \label{fig:spstorus}
\end{figure}

  We construct a diagram by using an arc  $a_n$ from each great circle $C_n$. For circle
  $C_n$, we use the arc that starts at $i_{n,n-p}$ and goes counterclockwise to
  $i_{n+p,n}$. The diagram for $T_{6,7}$ is shown in Figure
 ~\ref{fig:spstorusA}. 

  We can connect these arcs to form a closed loop by the same argument as in
  the previous proof, using the fact that $p$ and $q$ are coprime. Furthermore,
  for each $n$, the $2(p-1)$ points contained in the interior of arc $a_n$, 
  \[ i_{n,n-p+1}, \ldots, i_{n,n-1}, i_{n+1,n}, \ldots, i_{n+p-1,n}, \]
  are all intersections of this curve with itself. To see this, note that 
  an intersection of two circles $C_n$ and $C_m$ is equidistant from $n$ and $m$.
  Again, we get a diagram with $(p-1)q$ intersections.

  We choose the crossings for our diagram by letting arc $a_m$ be the overstrand
  for the $p-1$ crossings counterclockwise from $m$ and the understrand for the $p-1$
  crossings counterclockwise from $m'$ (see Figure
 ~\ref{fig:spstorusB}). We can then construct a conformation of $T_{p,q}$ that
  projects to our diagram, as in Figure~\ref{fig:spstorusC}.
 
  To show that $\sps[T_{q-1,q}]\geq q$ (and hence $\sps[T_{q-1,q}]=q$), we
  note that $\crn[T_{q-1,q}]=(q-2)q$, and apply the lower bound of Theorem
  \ref{thm:sps-crn}.
\end{proof}

\begin{rem}
  If the bound $\sbr[K]\le\sps[K]$ were to hold, the fact that $\sbr[T_{p,q}] =
  \min\{2p,q\}$ (see \cite{Kuiper}) would imply $\sps[T_{p,q}] = q$ for 
  $p < q < 2p$. 
\end{rem}

\subsection{Compositions of Trefoil Knots}
Another class of knots that we will examine is the set of compositions
of trefoil knots. Adams et al.~(see \cite{AdamsStickNumber}) computed the stick
index of such knots to be 
\[ \s[nT] = 2n+4. \]
We use this result to compute the planar stick index of compositions of 
trefoils.

\begin{proof}[Proof of Theorem \ref{thm:pltrefoil}] 
  By \cite{AdamsStickNumber} and Theorem \ref{thm:pl-s}, we know $\pl[nT] \le
  \s[nT]-1 = 2n+3$. Furthermore, it is known (see \cite{Schubert} or
  \cite{Schultens}) that the bridge index of a composition of knots is given
  by
  \[ \br[K_1\#K_2]=\br[K_1]+\br[K_2]-1. \]
  Since $\br[T]=2$, it follows that $\br[nT]=n+1$, and Theorem \ref{thm:pl-br}
  implies $\pl[nT] \ge 2\br[nT]+1 = 2n+3$.
\end{proof}

\begin{figure}
  \includegraphics[width=.25\textwidth]{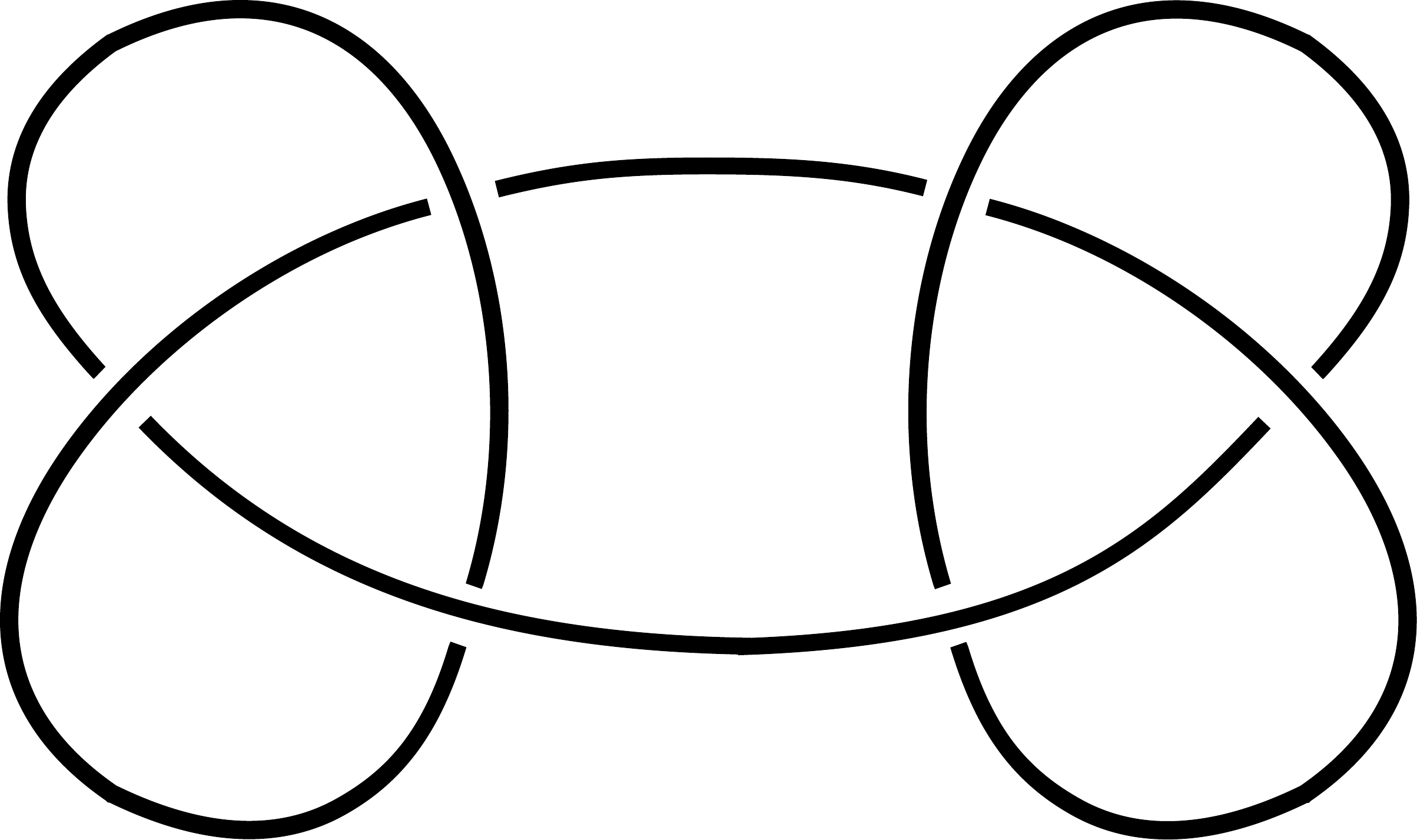}
  \hspace{.1\textwidth}
  \includegraphics[width=.25\textwidth]{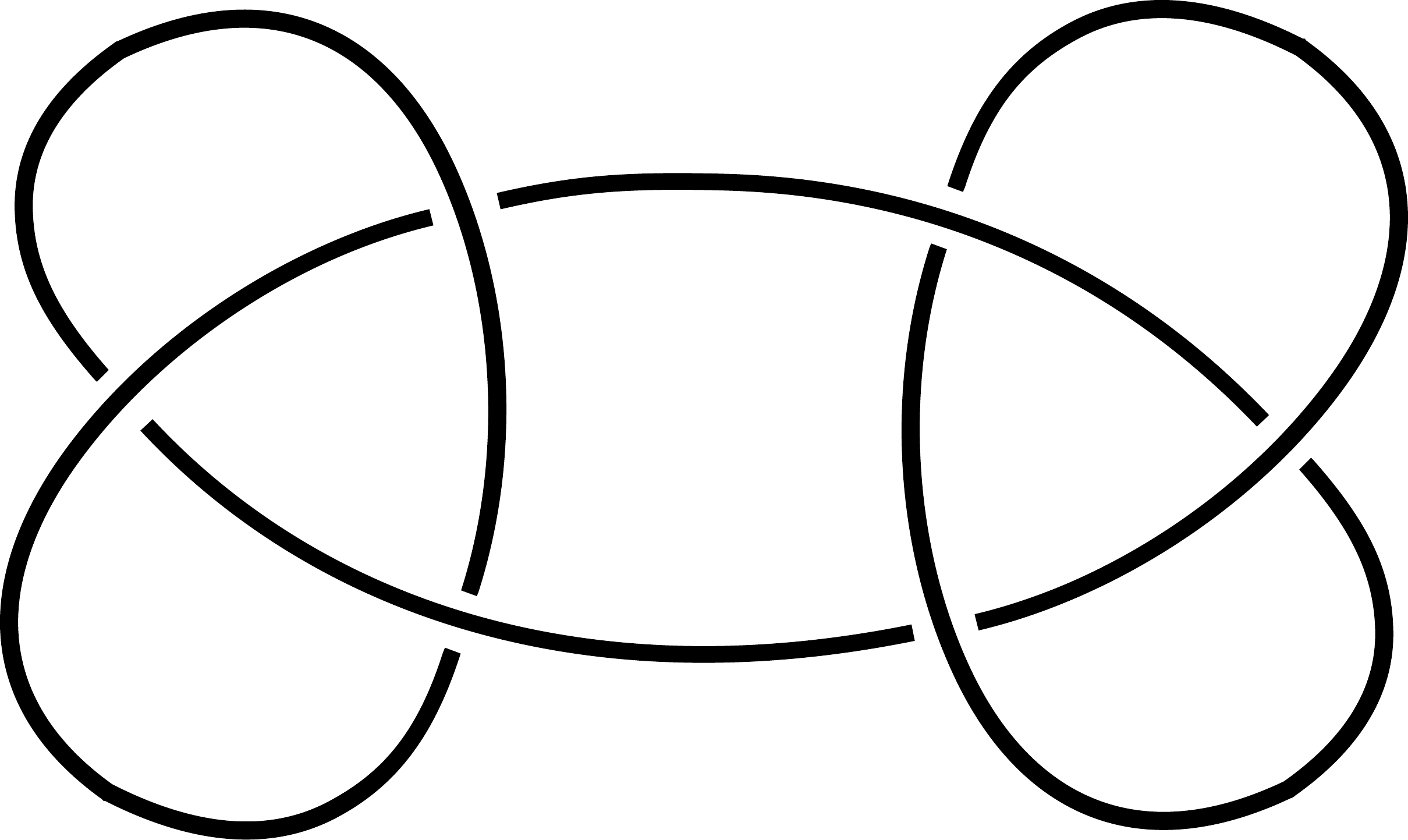}
  \caption{The square knot $T_L\# T_R$ and the granny knot $T_R\# T_R$.}
  \label{fig:handedness}
\end{figure}

When working with compositions of trefoils, we must be aware of a caveat.  The
trefoil knot is invertible, so there is a unique composition of two given
trefoil knots. However, since the trefoil is chiral, we must make a distinction
between left-handed and right-handed trefoil knots in compositions. For example, the
square and granny knots are the two distinct compositions of two trefoils (see
Figure~\ref{fig:handedness}). The $2n+4$-stick construction of $nT$ is
independent of handedness (as discussed in \cite{AdamsStickNumber}), so this
distinction does not affect Theorem \ref{thm:pltrefoil}.

The bulk of the proof of Theorem \ref{thm:spstrefoil} is in the following lemma:

\begin{lem} \label{lem:spstrefoil}
  If $m > 0$, we have 
  \begin{align*}
    \sps[mT_L\# mT_R] & \le 2m+2, \\
    \sps[(m+1)T_L\# mT_R] & \le 2m+3, \\
    \sps[mT_L\# (m+1)T_R] & \le 2m+3. 
  \end{align*}
\end{lem}
\begin{proof}
  We use the same conventions and notation as in the proof of Theorem
  \ref{thm:spstorus} and demonstrated in Figure~\ref{fig:spstorus_circles}. We
  start with $q$ great circles spaced symmetrically around the north pole, and
  label them from $C_1$ to $C_q$ counterclockwise. We assume throughout that $q > 2m+1$.
  We define the basepoint of
  circle $C_r$ as its farthest point from the origin, labelled $r$ and the point $r'$ as its
  closest point to the origin. Let $i_{r,s}$ denote the intersection of circles
  $C_r$ and $C_s$ that is counterclockwise from $r$ and clockwise from $s$.
  Again, consider all labels modulo $q$.

  Let $k$ be an integer, and $m = \lfloor k/2 \rfloor$. Then, we
  will prove the following set of statements for all $k$:
  
  \medskip

  \textit{For odd $k$, we have:} 
  \begin{enumerate}
    \item There is a diagram using $k+2$ arcs of circles that represents 
      $(m+1)T_L\# mT_R$.
    \item The diagram uses one arc from each circle corresponding to $-m,\ldots,m+2$,
      and we can pick an orientation on the diagram so that the circles are
      traversed in this order.
    \item If an arc of circle $C_r$ is in the diagram, it contains point 
      $r'$ but not the basepoint. 
    \item The vertices of the diagram are $i_{r,r+1}$ for 
      $p\in\{-m,\ldots,m+1\}$, along with $i_{m+2,-m}$. 
    \item Arc $a_{m+2}$ is the overstrand in all of its 
      crossings except for the crossing with arc $a_{-m}$.
  \end{enumerate}
  
  \medskip

  \textit{For even $k$, we have:} 
  \begin{enumerate}
    \item[($1'$)] There is a diagram using $k+2$ arcs of circles that 
      represents $mT_L\# mT_R$.
    \item[($2'$)] The diagram uses one arc from each circle corresponding to 
      $-m,\ldots,m+1$, and we can pick an orientation so that the circles are
      traversed in this order.
    \item[($3'$)] If an arc of circle $C_r$ is in the diagram, it contains 
      point $r'$ but not the basepoint. 
    \item[($4'$)] The vertices of the diagram are $i_{r,r+1}$ for 
      $p\in\{-m,\ldots,m\}$, along with $i_{m+1,-m}$. 
    \item[($5'$)] Arc $a_{-m}$ is the overstrand in all of its 
      crossings except for the crossing with arc $a_{m+1}$.
  \end{enumerate} 

  We start with the base case, $k=1$. We connect 
  the three vertices $i_{1,2}$, $i_{2,q}$, and $i_{2,1}$
  via arcs on circles corresponding to $1,2,q$ that pass through points $1',2',q'$
  respectively. This yields a three-crossing diagram, and by choosing
  crossings appropriately we obtain a left-handed trefoil that satisfies
  the conditions of the $k=1$ case of our induction hypothesis. The
  $q=7$ case is shown in Figure~\ref{fig:spstrefoil_induction}, and the picture
  looks similar for other $q$.

\begin{figure}
  \subfloat[]{\label{fig:spstrefoil_inductionA}
    \includegraphics[height=.27\textwidth]{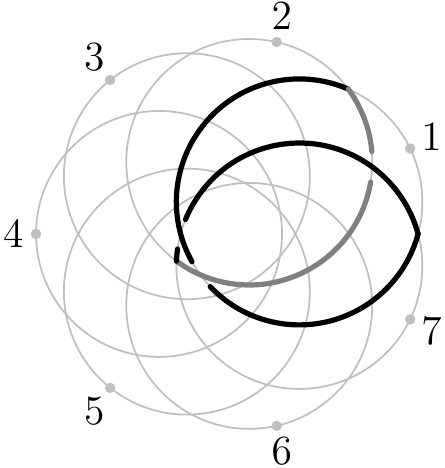}}
  \hspace{.05\textwidth}
  \subfloat[]{\label{fig:spstrefoil_inductionB}
    \includegraphics[height=.27\textwidth]{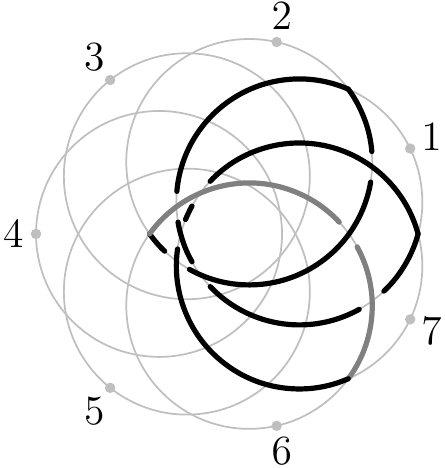}}
  \hspace{.05\textwidth}
  \subfloat[]{\label{fig:spstrefoil_inductionC}
    \includegraphics[height=.27\textwidth]{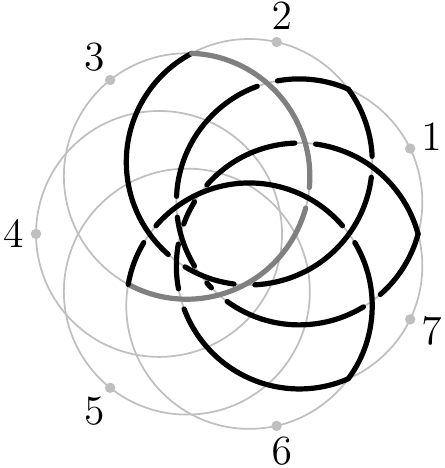}}
  \caption{Demonstrating the base case and the inductive step in the proof of 
    Lemma \ref{lem:spstrefoil}, with $n=7$. We go from the $k=1$ case to $k=2$
    to $k=3$.}
  \label{fig:spstrefoil_induction}
\end{figure} 

\begin{figure}[b]
  \includegraphics[height=.22\textwidth]{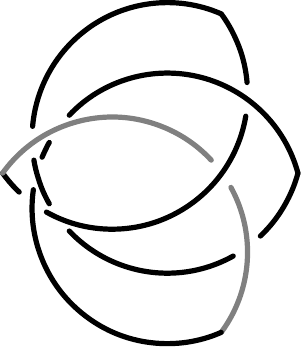}
  \hspace{.05\textwidth}
  \includegraphics[height=.22\textwidth]{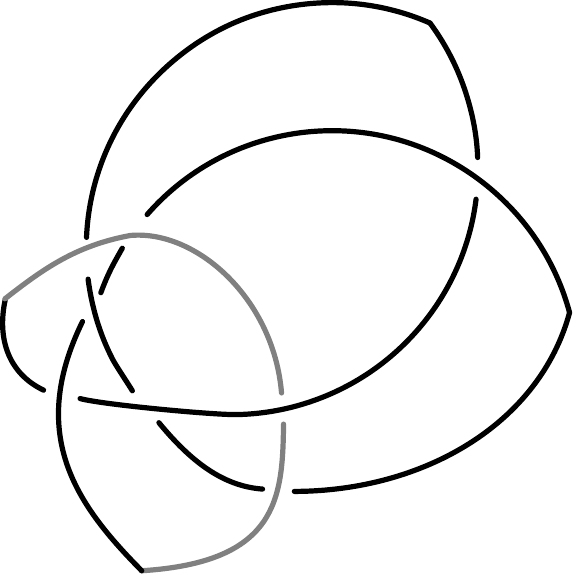}
  \hspace{.05\textwidth}
  \includegraphics[height=.20\textwidth]{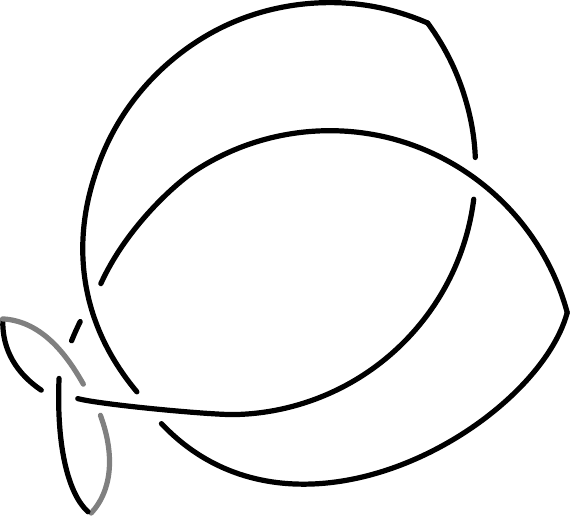}
  \caption{Demonstrating how the newly added arc can be modified via isotopy to 
    reveal a composed trefoil.}
  \label{fig:spstrefoil_isotopy}
\end{figure} 

  Suppose that the inductive hypotheses hold for some odd $k < q$. By
  assumption, we have a diagram of $(m+1)T_L\# mT_R$ using arcs corresponding to
  $-m,\ldots,m+2$. We modify our diagram by extending arc $a_{m+2}$ clockwise
  from $i_{m+1,m+2}$ to $i_{m+2,-m-1}$ and arc $a_{-m}$ clockwise from
  $i_{-m+1,-m}$ to $i_{-m-1,-m}$. We then add an arc of circle $C_{-m-1}$, which
  goes clockwise from $i_{m+2,-m-1}$ to $i_{-m-1,-m}$. Conditions (2) and (3)
  guarantee the new arcs of circles $C_{m+2}$ and $C_{-m}$ are extensions of the old
  ones.  Conditions ($2'$), ($3'$), and ($4'$) of the $k+1$ case immediately
  follow (see Figure~\ref{fig:spstrefoil_induction}).

  We must specify the crossings for our new diagram. We keep all crossings from
  the original diagram.  We let arc $a_{-m}$ be the overstrand at $i_{m+2,-m}$,
  arc $a_{m+2}$ be the overstrand at $i_{-m-1,m+2}$, and arc $a_{-m-1}$ be the
  overstrand at all of its other crossings (see Figure
 ~\ref{fig:spstrefoil_induction}).  By construction, condition ($5'$) is
  satisfied for the $k+1$ case.

  It remains to show that condition ($1'$) holds.  Because arc $a_{-m-1}$ is the
  overstrand at all crossings except its crossing with arc $a_{m+2}$, and arc $a_{m+2}$
  is the overstrand at all crossings except those with arc $a_{-m}$, we can move
  arc $a_{-m-2}$ as in Figure~\ref{fig:spstrefoil_isotopy}.  The result is the
  composition of a right-handed trefoil and the original diagram. By (1), this
  is $[(m+1)T_L\# mT_R]\# T_R$. This proves the inductive hypothesis for $k+1$.

  Finally, consider the case when $k$ is even. The argument is nearly identical
  to the previous one. This time, we extend arcs $a_{m+1}$ and $a_{-m-1}$ to the points
  $i_{m+1,m+2}$ and $i_{m+2,-m-1}$, and connect these points by adding an arc
  of circle $C_{m+2}$ between them.  We set arc $a_{m+2}$ as the overstrand in all of
  its crossings except that with arc $a_{-m-1}$, and make $a_{m+1}$ the overstrand at
  $i_{m+1,-m-1}$ (see Figure~\ref{fig:spstrefoil_induction}).  We can prove
  conditions (1)-(5) by the same arguments as above. This completes the 
  induction.

  Since $q$ is arbitrary, we have shown that $\sps[mT_L\# mT_R] \le 2m+2$ and \\
  $\sps[(m+1)T_L\# mT_R]\le 2m+3$ hold for all $m$.  Reflecting a $2m+3$-arc
  diagram of $mT_L\# (m+1)T_R$ yields a $2m+3$-arc diagram of $(m+1)T_L\#
  mT_R$, proving the last inequality.  
\end{proof}

\begin{figure}
  \subfloat[]{\label{fig:spstrefoil_vertexA}
    \includegraphics{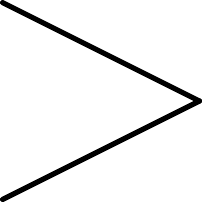}}
  \hspace{2cm}
  \subfloat[]{\label{fig:spstrefoil_vertexB}
    \includegraphics{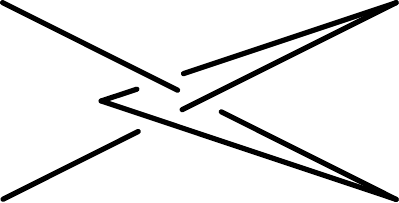}}
  \caption{Adding a trefoil at a vertex using two great circle arcs, which
    locally appear as line segments.}
  \label{fig:spstrefoil_vertex}
\end{figure}

\begin{proof}[Proof of Theorem \ref{thm:spstrefoil}]
  It remains to show that for $n > m$,
  \[ \sps[nT_L\# mT_R] =\sps[mT_L\# nT_R] \le 2n+1. \]
  We prove this by induction on $n$. The base case $n = m+1$ was proven
  in Lemma \ref{lem:spstrefoil}. For the inductive step, it suffices to show
  that given a knot, we can compose it with a trefoil by slightly extending two
  arcs past a vertex and adding two great circle arcs (see Figure
 ~\ref{fig:spstrefoil_vertex}).  
\end{proof}

\begin{rem}
  It would be interesting to know whether the bounds given in Theorem
  \ref{thm:spstrefoil} are tight.  If $\sbr[K]\le\sps[K]$ were to hold, then
  the bounds of Lemma \ref{lem:spstrefoil} would be tight, but this would not
  help with the general case.
\end{rem}

%%%%%%%%%%%%%%%%%%%%%%%%%%%%%%%%%%%%%%%%%%%%%%%%%%%%%%%%%%%%%%%%%%%%%%%%%%%%%%
%%%%%%%%%%%%%%%%%%%%%%%%%%%%%%%%%%%%%%%%%%%%%%%%%%%%%%%%%%%%%%%%%%%%%%%%%%%%%%
\section{Classification of Knots with $\sps=4$.} \label{section:sps4}

\begin{figure}
  \includegraphics[height=.2\textwidth]{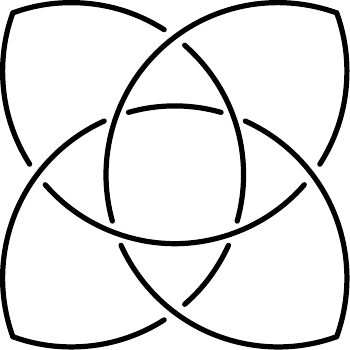}
  \hspace{.2\textwidth}
  \includegraphics[height=.2\textwidth]{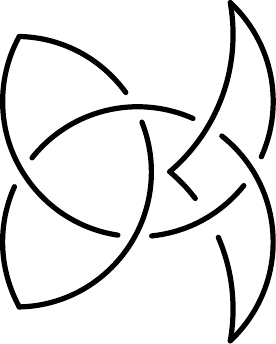}
  \caption{Spherical projective stick diagrams of the square knot $T_L\# T_R$
    and the granny knot $T_L\# T_L$ that realize their spherical stick indices.}
  \label{fig:spssqgran}
\end{figure}

Our construction in Theorem \ref{thm:spstrefoil} depends on the handedness of
the composed trefoils, but this does not imply that spherical stick index
depends on handedness. However, we will show that $T_L\# T_R$ and $T_L\# T_L$
have different spherical stick indices (see Figure~\ref{fig:spssqgran}). 
By Theorem \ref{thm:sps-crn}, any knot with at
least 4 crossings has $\sps[K]\ge 4$. Combining this with Theorem 
\ref{thm:spstrefoil} gives $\sps[T_L\# T_R] = 4$ and 
$\sps[T_L\# T_L]\le 5$. To show that $\sps[T_L\# T_L]=5$, we classify all knots
with $\sps[K]\le 4$. 

\begin{figure}
  \includegraphics[width=.25\textwidth]{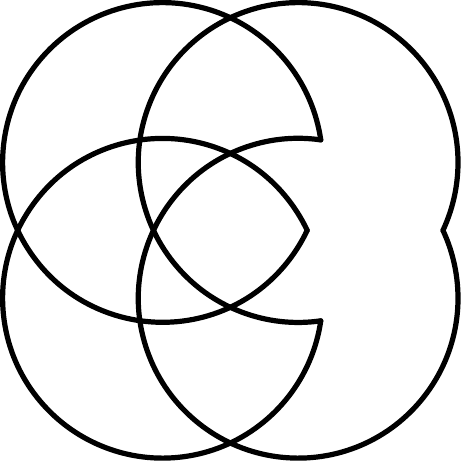} 
  \hspace{.05\textwidth}
  \includegraphics[width=.25\textwidth]{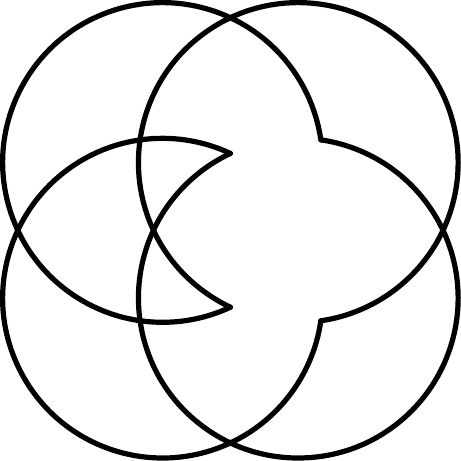} 
  \hspace{.05\textwidth}
  \includegraphics[width=.25\textwidth]{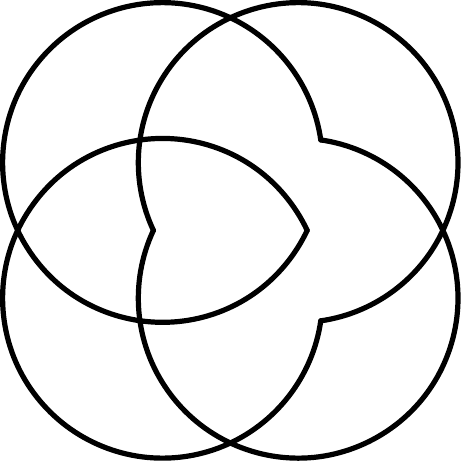} 
  \caption{The three possible loops with at least $6$ crossings.}
  \label{fig:sps4_loops}
\end{figure}

\begin{proof}[Proof of Theorem \ref{thm:sps4classification}]
To construct an $\sps$-4 diagram, we start with a configuration of four great
circles on the sphere, as in Figure~\ref{fig:sps4_circles}.  It is not
difficult to show that any generic configuration divides the sphere into
triangular and quadrilateral regions, and that no two triangles or
quadrilaterals can share a side. The only way to do this on a sphere (up to
isotopy) is the arrangement shown in Figure~\ref{fig:sps4_circles}. Thus, this
is the only configuration that we need to consider. 

\begin{figure}[b]
  \includegraphics[width=.3\textwidth]{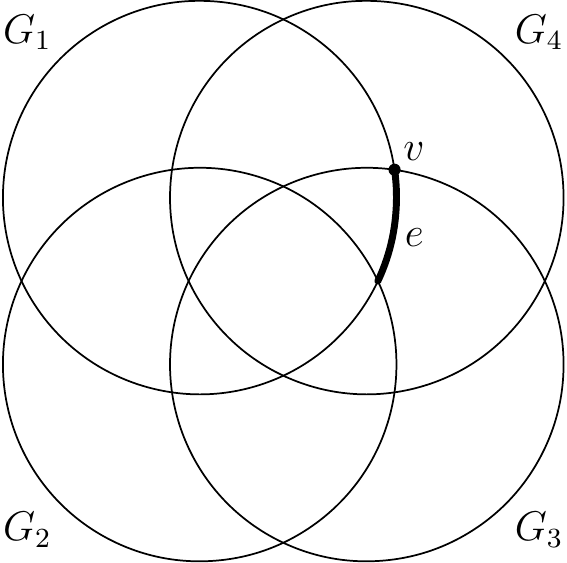} 
  \caption{The stereographic projection of an arrangement of four great circles 
  on a sphere.} 
  \label{fig:sps4_circles}
\end{figure}

Given this diagram, we find all ways to form a closed loop out of one arc from
each great circle $G_i$. We note that the complement of such a closed loop,
which consists of all arcs \emph{removed} from the diagram to obtain the first
loop, is itself a closed loop with one arc from each great circle. It will be
convenient to describe loops by their complements, because the loops with the
most crossings have simple complements.

We choose a closed curve in the diagram. Suppose the complement contains $n$
edges, and hence $n$ vertices. Four of these vertices are the ``turning
vertices'' where the complement switches from one circle to another, and the
other $n-4$ are ``passing vertices'' where the complement passes straight
through an intersection.  Note that the complement may pass through a single
vertex twice, removing all four edges. Thus, at least $\lceil(n-4)/2\rceil$
vertices are passed through. These vertices, along with the 4 turning vertices,
are not crossings of the original curve, so the number of remaining crossings
is at most
\[ 12 - ( 4 + \lceil (n-4)/2 \rceil) = 10 - \lceil n/2 \rceil. \]

Note that for $n > 8$, there will be at most five crossings. We will see that we pick up all knots of five or fewer crossings from the complements with $n \le 6$.
 
By symmetry, it is not hard to check that all pairs of an edge $e$ and an
adjacent vertex $v$ are combinatorially equivalent. Thus, we need only consider
complementary loops including $v$ as a turning vertex and containing $e$, as in Figure \ref{fig:sps4_circles}.

By explicitly considering complements constructed from 7 or 8 edges, we can see that they only produce knots of four or fewer crossings. Hence, we can limit consideration to complements of 4, 5, or 6 edges. There are none with 5 edges as any such would need to have two edges on one great circle and one on each of the others, and we cannot close such a complement up. Up to equivalence of diagrams we obtain one complement with four edges and two with six edges, as appear in 
Figure~\ref{fig:sps4_loops}. 

We consider all possible crossing choices for the diagrams in Figure 
\ref{fig:sps4_loops} and identify what knots result. This list of knots 
includes all knots of five or fewer crossings. We conclude that the nontrivial knots with 
$\sps[K]\le 4$ are 
\[ 3_1,4_1,5_1,5_2,6_1,6_2,6_3,7_4,8_{18},8_{19},8_{20}, T_L\# T_R. \]
It follows from Theorem \ref{thm:sps-crn} that the trefoil has $\sps[K]=3$
and that the other knots listed have $\sps[K]=4$.
\end{proof}

\begin{table}
\begin{tabular}{|c|c|c|c|}
  \hline \rule{0pt}{.205\textwidth} 
  \includegraphics[height=.2\textwidth]{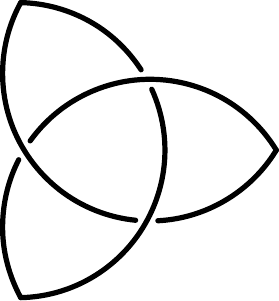} 
  & \includegraphics[width=.2\textwidth]{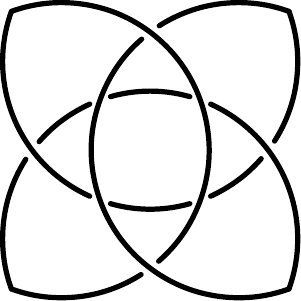} 
  & \includegraphics[width=.2\textwidth]{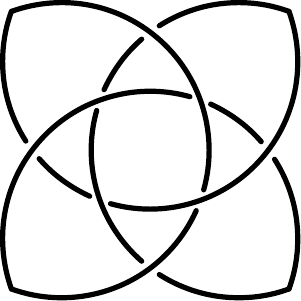} 
  & \includegraphics[width=.2\textwidth]{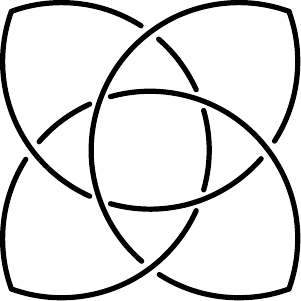} 
  \\ $3_1$ & $4_1$ & $5_1$ & $5_2$
  \\ \hline \rule{0pt}{.205\textwidth} 
  \includegraphics[width=.2\textwidth]{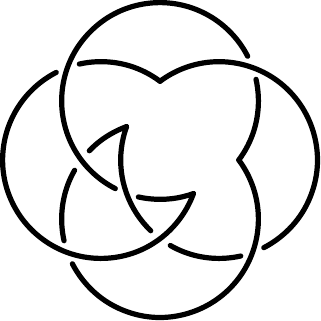} 
  & \includegraphics[width=.2\textwidth]{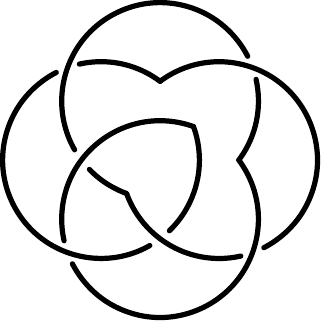} 
  & \includegraphics[width=.2\textwidth]{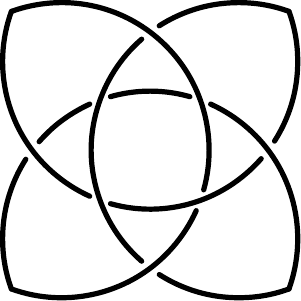} 
  & \includegraphics[width=.2\textwidth]{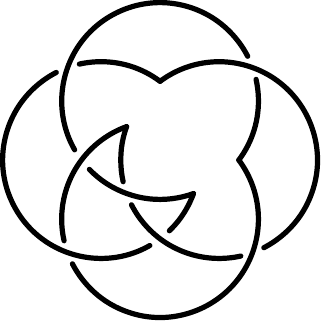} 
  \\ $6_1$ & $6_2$ & $6_3$ & $7_4$ 
  \\ \hline \rule{0pt}{.205\textwidth} 
  \includegraphics[width=.2\textwidth]{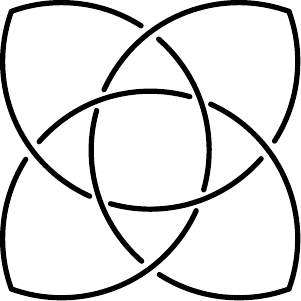} 
  & \includegraphics[width=.2\textwidth]{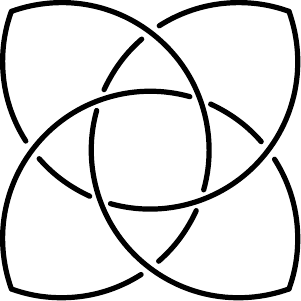} 
  & \includegraphics[width=.2\textwidth]{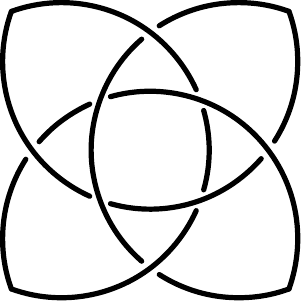} 
  & \includegraphics[width=.2\textwidth]{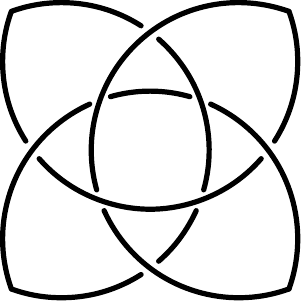} 
  \\ $8_{18}$ & $8_{19}$ & $8_{20}$ & Square Knot \\
  \hline
\end{tabular}
\caption{Nontrivial knots with $\sps\le 4$.}
\label{fig:sps4table}
\end{table}

Using a computer (with the help of the program Knotscape), we carried
out a similar process to classify knots with $\sps[K] = 5$. We found
that there are $666$ prime knots and $17$ composite knots with $\sps[K]=5$. 
In particular, we found that all knots of eight or fewer
crossings (prime or composite) have $\sps[K]\le 5$. The list also includes
nine-crossing knots except for $9_2$, $9_3$, $9_4$,
$9_{15}$, $9_{18}$, $9_{23}$, $9_{36}$, $4_1\# 5_1$, $4_1\# 5_2$, and
$T_L\# T_L\# T_L$ and all ten-crossing nonalternating prime knots except for 
$10_{152}$ and $10_{154}$. 

Since $\sps[T_L\# T_L\# T_L] > 5$, spherical stick index distinguishes between
the two distinct types of compositions of three trefoils, one consisiting of all left or all right trefoils and one with a mixture of the two. We found a few torus knots
with $\sps[T_{p,q}]$ strictly less than $q$: $\sps[T_{2,5}] = 4$ and
$\sps[T_{2,7}] = \sps[T_{2,9}] = 5$. Note that these knots still satisfy
$\sbr[K]\le\sps[K]$, as $\sbr[T_{p,q}] = \min\{ 2p,q\}$ is $4$ in these cases.

\begin{proof}[Proof of Corollary \ref{cor:spscomposition}]
  From above, $\sps[T_L\# T_L] = \sps[T_L\# T_L\# T_R] = 5$.
\end{proof}

\bibliographystyle{alpha}
\bibliography{spsbib}

\end{document}